\documentclass[11pt]{article}
\usepackage{latexsym,amssymb,amsmath,geometry,cite,enumerate,float}
\newtheorem{theorem}{Theorem}
\newtheorem{lemma}[theorem]{Lemma}
\newtheorem{proposition}[theorem]{Proposition}

\usepackage{setspace}

\begin{document}

%\linenumbers
\onehalfspace

\title{Degenerate Matchings and Edge Colorings}

\author{Julien Baste$^1$ and Dieter Rautenbach$^2$}

\date{}

\maketitle

\begin{center}
$^1$ LIRMM, Universit\'{e} de Montpellier, Montpellier, France,
\texttt{baste@lirmm.fr}\\[3mm]
$^2$ Institute of Optimization and Operations Research, Ulm University, Ulm, Germany, \texttt{dieter.rautenbach@uni-ulm.de}
\end{center}

\begin{abstract}
A matching $M$ in a graph $G$ is $r$-degenerate if the subgraph of $G$ induced by the set of vertices incident with an edge in $M$ is $r$-degenerate. 
Goddard, Hedetniemi, Hedetniemi, and Laskar 
(Generalized subgraph-restricted matchings in graphs, Discrete Mathematics 293 (2005) 129-138) introduced the notion of acyclic matchings, which coincide with $1$-degenerate matchings.
Solving a problem they posed, we describe an efficient algorithm to determine the maximum size of an $r$-degenerate matching in a given chordal graph. Furthermore, we study the $r$-chromatic index of a graph defined as the minimum number of $r$-degenerate matchings into which its edge set can be partitioned,
obtaining upper bounds and discussing extremal graphs.
\end{abstract}

{\small 

\begin{tabular}{lp{13cm}}
{\bf Keywords:} & Matching; edge coloring; induced matching; acyclic matching; uniquely restricted matching
\end{tabular}
}

\section{Introduction}

Matchings in graphs are a central topic of graph theory and combinatorial optimization \cite{lp}.
While classical matchings are tractable, several well known types of more restricted matchings,
such as induced matchings \cite{c,sv}
or uniquely restricted matchings \cite{ghl},
lead to hard problems.
Goddard, Hedetniemi, Hedetniemi, and Laskar \cite{ghhl} proposed to study so-called subgraph-restricted matchings in general.
In particular, they introduce the notion of acyclic matchings.
By a simple yet elegant argument (cf. Theorem 4 in \cite{ghhl}) 
they show that finding a maximum acyclic matching in a given graph is hard in general,
and they explicitly pose the problem to describe a fast algorithm for the acyclic matching number in interval graphs.
In the present paper, we solve this problem for the more general chordal graphs.
Furthermore, we study the edge coloring notion corresponding to acyclic matchings.

Before we give exact definitions and discuss our results as well as related research, 
we introduce some terminology.
We consider finite, simple, and undirected graphs, and use standard notation.
A {\it matching} in a graph $G$ is a subset $M$ of the edge set $E(G)$ of $G$ 
such that no two edges in $M$ are adjacent.
Let $V(M)$ be the set of vertices incident with an edge in $M$.
$M$ is {\it induced} \cite{c} 
if the subgraph $G[V(M)]$ of $G$ induced by the set $V(M)$ is $1$-regular,
that is, $M$ is the edge set of $G[V(M)]$.
Induced matching are also known as {\it strong} matchings.
$M$ is {\it uniquely restricted} \cite{ghl} 
if there is no other matching $M'$ in $G$ distinct from $M$ 
that satisfies $V(M)=V(M')$.
It is easy to see that $M$ is uniquely restricted if and only if there is no $M$-alternating cycle in $G$,
which is a cycle in $G$ every second edge of which belongs to $M$ \cite{ghl}.
Finally, $M$ is {\it acyclic} \cite{ghhl} if $G[V(M)]$ is a forest.
Let 
$\nu(G)$, 
$\nu_s(G)$,
$\nu_{ur}(G)$, and 
$\nu_1(G)$
be the maximum sizes of 
a matching,
an induced matching,
a uniquely restricted matching, and 
an acyclic matching in $G$, respectively.
Since every induced matching is acyclic, 
and every acyclic matching is uniquely restricted, we have
$$\nu_s(G)\leq \nu_1(G)\leq \nu_{ur}(G)\leq \nu(G).$$
We chose the notation ``$\nu_1(G)$'' rather than something like ``$\nu_{ac}(G)$'', 
because we consider some further natural generalization.

For a non-negative integer $r$, a graph $G$ is {\it $r$-degenerate} 
if every subgraph of $G$ of order at least one has a vertex of degree at most $r$. 
Note that a graph is a forest if and only if it is $1$-degenerate.
An {\it $r$-degenerate order} of a graph $G$ is a linear order $u_1,\ldots,u_n$ of its vertices 
such that, for every $i$ in $[n]$, the vertex $v_i$ has degree at most $r$ in $G[\{ v_i,\ldots,v_n\}]$,
where $[n]$ is the set of the positive integers at most $n$.
Clearly, a graph is $r$-degenerate if and only if it has an $r$-degenerate order.

Now, let a matching $M$ in a graph $G$ be {\it $r$-degenerate} 
if the induced subgraph $G[V(M)]$ is $r$-degenerate,
and let $\nu_r(G)$ denote the maximum size of an $r$-degenerate matching in $G$.

\medskip

\noindent For every type of matching, there is a corresponding edge coloring notion.
An {\it edge coloring} of a graph $G$ is a partition of its edge set into matchings.
An edge coloring is 
{\it induced (strong)},
{\it uniquely restricted}, and 
{\it $r$-degenerate}
if each matching in the partition has this property, respectively.
Let 
$\chi'(G)$,
$\chi'_s(G)$,
$\chi'_{ur}(G)$, and
$\chi'_r(G)$ 
be the minimum numbers of colors needed for the corresponding colorings, respectively.
Clearly, 
$$\chi'_s(G)\geq \chi'_1(G)\geq \chi'_{ur}(G)\geq \chi'(G).$$
In view of the hardness of the restricted matching notions,
lower bounds on the matching numbers \cite{hr,hqt,jrs,kmm},
upper bounds on the chromatic indices \cite{a,brs},
efficient algorithms for restricted graph classes \cite{c2,lo,fjj,cw,djprs}, and 
approximation algorithms have been studied \cite{brs,r}.
There is only few research concerning acyclic matchings;
Panda and Pradhan \cite{pp} describe efficient algorithms for chain graphs and bipartite permutation graphs. 

Vizing's \cite{v} famous theorem says that the {\it chromatic index} $\chi'(G)$ of $G$
is either $\Delta(G)$ or $\Delta(G)+1$, 
where $\Delta(G)$ is the maximum degree of $G$.
Induced edge colorings have attracted much attention 
because of the conjecture $\chi_s'(G)\leq \frac{5}{4}\Delta(G)^2$
posed by Erd\H{o}s and Ne\v{s}et\v{r}il (cf. \cite{fsgt}).
Building on earlier work of Molloy and Reed \cite{mr},
Bruhn and Joos \cite{bj} showed $\chi_s'(G)\leq 1.93\Delta(G)^2$ provided that $\Delta(G)$ is sufficiently large. 
In \cite{brs} we showed $\chi'_{ur}(G)\leq \Delta(G)^2$ 
with equality if and only if $G$ is the complete bipartite graph $K_{\Delta(G),\Delta(G)}$.

\medskip

\noindent Our results are upper bounds on $\chi'_r(G)$ with the discussion of extremal graphs, and an efficient algorithm for $\nu_r(G)$ in chordal graphs,
solving the problem posed in \cite{ghhl}.

\section{Bounds on the $r$-degenerate chromatic index}

Since, for every two positive integers $r$ and $\Delta$,
every $r$-degenerate matching of the complete bipartite graph $K_{\Delta,\Delta}$ of order $2\Delta$ has size at most $r$, we obtain $\chi_r'(K_{\Delta,\Delta})\geq \frac{\Delta^2}{r}$.

Our first result gives an upper bound in terms of $r$ and $\Delta$.

\begin{theorem}\label{theorem1}
If $r$ is a positive integer and $G$ is a graph of maximum degree at most $\Delta$, then
\begin{eqnarray}\label{e1}
\chi_r'(G)\leq \frac{2(\Delta-1)^2}{r+1}+2(\Delta-1)+1.
\end{eqnarray}
\end{theorem}
{\it Proof:}
Let $K=\left\lfloor\frac{2(\Delta-1)^2}{r+1}+2(\Delta-1)+1\right\rfloor$.
The proof is based on a inductive coloring argument. 
We may assume that all but exactly one edge $uv$ of $G$ are colored using colors in $[K]$
such that, for every color $\alpha$ in $[K]$,
the edges of $G$ colored with $\alpha$ form an $r$-degenerate matching.
We consider the colors in $[K]$ that are forbidden by colors of the edges close to $uv$.
In order to complete the proof, 
we need to argue that there is always still some available color for $uv$ in $[K]$.

\begin{figure}[H]
\begin{center}
%TeXCAD Picture [1.pic]. Options:
%\grade{\on}
%\emlines{\off}
%\epic{\off}
%\beziermacro{\on}
%\reduce{\on}
%\snapping{\on}
%\pvinsert{% Your \input, \def, etc. here}
%\quality{8.000}
%\graddiff{0.005}
%\snapasp{1}
%\zoom{11.3137}
\unitlength 2mm % = 2.845pt
\linethickness{0.4pt}
\ifx\plotpoint\undefined\newsavebox{\plotpoint}\fi % GNUPLOT compatibility
\begin{picture}(63,31)(0,0)
\put(26,18){\circle*{1}}
\put(46,18){\circle*{1}}
\put(36,6){\oval(20,4)[]}
\put(26,18){\line(-1,0){5}}
\put(26,18){\line(1,0){25}}
\put(46,18){\line(-1,-2){5}}
\put(26,18){\line(1,-2){5}}
\put(19,18){\oval(4,20)[]}
\put(53,18){\oval(4,20)[]}
\put(19,26){\circle*{1}}
\put(19,10){\circle*{1}}
\put(53,10){\circle*{1}}
\put(26,20){\makebox(0,0)[cc]{$u$}}
\put(46,20){\makebox(0,0)[cc]{$v$}}
\put(19,18){\makebox(0,0)[cc]{$N_u$}}
\put(53,18){\makebox(0,0)[cc]{$N_v$}}
\put(36,6){\makebox(0,0)[cc]{$N_{u,v}$}}
\put(19,13){\makebox(0,0)[cc]{$u'$}}
\put(53,13){\makebox(0,0)[cc]{$v'$}}
\put(19,23){\makebox(0,0)[cc]{$w$}}
\put(21,24){\line(5,-6){5}}
\put(26,18){\line(-5,-6){5}}
\put(51,24){\line(-5,-6){5}}
\put(46,18){\line(5,-6){5}}
\put(26,18){\line(3,-2){15}}
\put(46,18){\line(-3,-2){15}}
\put(24,24){\makebox(0,0)[cc]{$E_u$}}
\put(48,24){\makebox(0,0)[cc]{$E_v$}}
\put(36,15){\makebox(0,0)[cc]{$E_{u,v}$}}
\put(19,26){\line(-1,0){10}}
\put(10,5){\line(-1,0){10}}
\put(9,21){\line(2,1){10}}
\put(0,0){\line(2,1){10}}
\put(19,26){\line(-2,1){10}}
\put(10,5){\line(-2,1){10}}
\put(14,31){\makebox(0,0)[cc]{$E_w$}}
\put(5,10){\makebox(0,0)[cc]{$E_{u''}$}}
\put(53,10){\line(1,0){10}}
\put(63,15){\line(-2,-1){10}}
\put(53,10){\line(2,-1){10}}
\put(59,16){\makebox(0,0)[cc]{$E_{v'}$}}
\put(10,10){\oval(4,14)[]}
\put(10,14){\makebox(0,0)[cc]{$N_{u'}$}}
\put(10,5){\circle*{1}}
%\emline(12,14)(19,10)
\multiput(12,14)(.058823529,-.033613445){119}{\line(1,0){.058823529}}
%\end
%\emline(19,10)(12,6)
\multiput(19,10)(-.058823529,-.033613445){119}{\line(-1,0){.058823529}}
%\end
\put(12,10){\line(1,0){7}}
\put(10,8){\makebox(0,0)[cc]{$u''$}}
\put(28,6){\circle*{1}}
\put(44,6){\circle*{1}}
\put(53,26){\circle*{1}}
\end{picture}
\end{center}
\caption{Vertices and edges close to $uv$. The indicated objects $N_{u'}$, $u''$, and $E_{u''}$ will be introduced and discussed in the proof of Theorem \ref{theorem2}. Note that the set $N_u\cup N_v\cup N_{u,v}$ is not required to be independent,
that is, the sets $E_w$ and $E_{w'}$ may intersect 
for distinct vertices $w$ and $w'$ in $N_u\cup N_v\cup N_{u,v}$.}\label{fig1}
\end{figure}
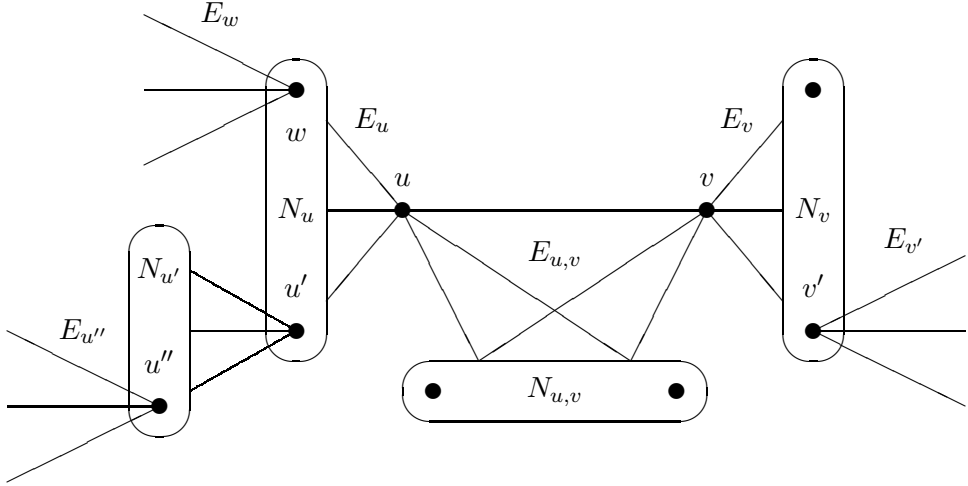

We introduce some notation illustrated in Figure \ref{fig1}.
Let 
$N_u=N_G(u)\setminus N_G[v]$,
$N_v=N_G(v)\setminus N_G[u]$, and
$N_{u,v}=N_G(u)\cap N_G(v)$.
Let
$n_u=|N_u|$,
$n_v=|N_v|$, and
$n_{u,v}=|N_{u,v}|$.
Clearly, 
$n_u+n_{u,v}=d_G(u)-1\leq \Delta-1$
and
$n_v+n_{u,v}=d_G(v)-1\leq \Delta-1$.
Let 
$E_u$ be the set of edges between $u$ and $N_u$,
$E_v$ be the set of edges between $v$ and $N_v$,
$E_{u,v}$ be the set of edges between $\{ u,v\}$ and $N_{u,v}$, and,
for every vertex $w\in N_u\cup N_v\cup N_{u,v}$, 
let $E_w$ be the set of edges incident with $w$ but not incident with $u$ or $v$.
Clearly, 
$|E_u|+|E_v|+|E_{u,v}|=(d_G(u)-1)+(d_G(v)-1)\leq 2(\Delta-1)$
and 
$|E_w|\leq \Delta-1$ for every vertex $w\in N_u\cup N_v\cup N_{u,v}$.

\noindent Let $F_1$ be the colors that appear on edges in $E_u\cup E_v\cup E_{u,v}$.
Clearly, every color in $F_1$ is forbidden for $uv$, 
because each color class must be a matching.
Let $F_2$ be the colors $\alpha$ in $[K]$ that do not belong to $F_1$
such that 
$$d_u^{\alpha}+2d_{u,v}^{\alpha}+d_v^{\alpha}\geq r+1,$$
where 
$d_u^{\alpha}$ is the number of vertices in $N_u$ incident with an edge colored $\alpha$,
$d_v^{\alpha}$ is the number of vertices in $N_v$ incident with an edge colored $\alpha$, and
$d_{u,v}^{\alpha}$ is the number of vertices in $N_{u,v}$ incident with an edge colored $\alpha$.
Note that, since $F_1$ and $F_2$ are disjoint, 
none of the edges contributing to $d_u^{\alpha}+2d_{u,v}^{\alpha}+d_v^{\alpha}$
is incident with $u$ or $v$.

If there is some $\alpha$ in $[K]\setminus (F_1\cup F_2)$, 
then neither $u$ nor $v$ is incident with an edge of color $\alpha$,
and $d_u^{\alpha}+2d_{u,v}^{\alpha}+d_v^{\alpha}\leq r$.
This implies
$\min\{ d_v^{\alpha}+d_{u,v}^{\alpha},d_v^{\alpha}+d_{u,v}^{\alpha}\}\leq \lfloor r/2\rfloor\leq r-1$
and 
$\max\{ d_v^{\alpha}+d_{u,v}^{\alpha},d_v^{\alpha}+d_{u,v}^{\alpha}\}\leq r$.
Hence, coloring $uv$ with color $\alpha$, 
the edges of $G$ colored $\alpha$ form an $r$-degenerate matching.
As explained above this would complete the proof.
Therefore, we may assume that $F_1\cup F_2=[K]$.

Note that 
\begin{eqnarray}
|F_1| & \leq & |E_u\cup E_v\cup E_{u,v}|\label{e3a}\\
&=&(d_G(u)-1)+(d_G(v)-1)\nonumber\\
&\leq &2(\Delta-1)\label{e3b}
\end{eqnarray}
with equality if and only if 
\begin{enumerate}
\item[(a)] all edges in $E_u\cup E_v\cup E_{u,v}$ are colored differently (equality in (\ref{e3a})), and 
\item[(b)] $u$ and $v$ have degree $\Delta$ (equality in (\ref{e3b})).
\end{enumerate}
Furthermore,
\begin{eqnarray}
(r+1)|F_2| & \leq & \sum\limits_{\alpha\in F_2}\left(d_u^{\alpha}+2d_{u,v}^{\alpha}+d_v^{\alpha}\right)\label{e2a}\\
& \leq & \sum\limits_{w\in N_u}|E_w|+2\sum\limits_{w\in N_{u,v}}|E_w|+\sum\limits_{w\in N_v}|E_w|\label{e2b}\\
& \leq & (\Delta-1)n_u+2(\Delta-2)n_{u,v}+(\Delta-1)n_v\label{e2c}\\
& \leq & (\Delta-1)(d_G(u)-1)+(\Delta-1)(d_G(v)-1)\label{e2d}\\
& \leq & 2(\Delta-1)^2.\label{e2e} 
\end{eqnarray}
Note that $(r+1)|F_2|=2(\Delta-1)^2$ if and only if equality holds in (\ref{e2a}) to (\ref{e2e}),
which implies that 
\begin{enumerate}
\item[(c)] $d_u^{\alpha}+2d_{u,v}^{\alpha}+d_v^{\alpha}=r+1$ for every color $\alpha$ in $F_2$ 
(equality in (\ref{e2a})),
\item[(d)] all edges in $\bigcup\limits_{w\in N_u\cup N_{u,v}\cup N_w}E_w$ have a color from $F_2$
(equality in (\ref{e2b})),
\item[(e)] all vertices in $N_u\cup N_{u,v}\cup N_v$ have degree $\Delta$ (equality in (\ref{e2c})),
\item[(f)] $n_{u,v}=0$, that is, $u$ and $v$ have no common neighbor (equality in (\ref{e2d})), and
\item[(g)] $u$ and $v$ have degree $\Delta$ (equality in (\ref{e2e})).
\end{enumerate}
Altogether, we obtain
\begin{eqnarray*}
|K| & = & |F_1\cup F_2| = |F_1|+|F_2|\leq 2(\Delta-1)+\frac{2(\Delta-1)^2}{\alpha+1},   
\end{eqnarray*}
contradicting the choice of $K$.
This completes the proof. $\Box$

\medskip

\noindent For $r=1$, the bound from Theorem \ref{theorem1} simplifies to $\Delta^2$.
In view of $K_{\Delta,\Delta}$, Theorem \ref{theorem1} is tight in this case, and, as we show next, $K_{\Delta,\Delta}$ is the only extremal graph.

\begin{theorem}\label{theorem2}
If $G$ is a graph of maximum degree at most $\Delta$,
then $\chi_1'(G)=\Delta^2$ if and only if $G$ is $K_{\Delta,\Delta}$.
\end{theorem}
{\it Proof:} By Theorem \ref{theorem1}, we have $\chi_1'(G)\leq \Delta^2$.
Since $\chi_1'(K_{\Delta,\Delta})=\Delta^2$,
it suffices to show that $\chi_1'(G)=\Delta^2$ implies that $G$ is $K_{\Delta,\Delta}$.
Therefore, we consider a $1$-degenerate edge coloring of $G$ using colors in $[\Delta^2]$
such that the number of edges colored $\Delta^2$ is as small as possible.
Let $uv$ be an edge colored $\Delta^2$.

We use the notation and observations 
from the proof of Theorem \ref{theorem1}.
Recall that 
$|F_1|\leq 2(\Delta-1)$ and $2|F_2|\leq 2(\Delta-1)^2$,
which implies
$|F_1\cup F_2|\leq \Delta^2-1$.
Furthermore, 
recall that $uv$ can be colored with any color in $[\Delta^2-1]\setminus (F_1\cup F_2)$.
By the choice of the coloring, 
these observations imply that 
$F_1\cup F_2=[\Delta^2-1]$,
$|F_1|=2(\Delta-1)$, and 
$2|F_2|=2(\Delta-1)^2$.
The latter two equalities imply that the properties (a) to (f) hold.
%In particular, by (c), for every color $\alpha$ in $F_2$,
%we have $d_u^{\alpha}+d_v^{\alpha}=2$,
%that is, exactly two vertices in $N_u\cup N_v$ are incident with an edge colored %$\alpha$.

\medskip

\noindent Let $u'\in N_u$ and let $\alpha$ be the color of the edge $uu'$.
We introduce some more notation already illustrated in Figure \ref{fig1}.
Let $N_{u'}=N_G(u')\setminus \{ u\}$.
For every vertex $w$ in $N_{u'}$, 
let $E_w$ be the set of edges incident with $w$ but not incident with $u'$.
Let $E^2_{u'}=\bigcup\limits_{w\in N_{u'}}E_w$.

For every color $\beta$ in $[\Delta^2-1]$,
let $k_{\beta}$ be the number of vertices $w$ in $\{ v\}\cup N_v$ 
such that $E_w$ contains an edge colored $\beta$,
and, similarly, 
let $k_{\beta}'$ be the number of vertices $w$ in $\{ u'\}\cup N_{u'}$ 
such that $E_w$ contains an edge colored $\beta$.
Since the color classes are matchings, for every such color $\beta$, 
each of the sets $E_w$ for $w$ in $\{ v,u'\}\cup N_v\cup N_{u'}$
contains at most one edge colored $\beta$.
By (a), (c), and (d),
all edges in $E_v$ have different colors from $F_1$,
all edges in $\bigcup\limits_{w\in N_v}E_w$ have colors from $F_2$,
$k_{\beta}\in \{ 0,1\}$ for every color $\beta$ in $F_1$, and
$k_{\beta}\in \{ 0,1,2\}$ for every color $\beta$ in $F_2$.
By (b), (e), and (f), 
we obtain 
$\sum\limits_{\beta\in [\Delta^2-1]}k_{\beta}
=|E_v|+\sum\limits_{w\in N_v}|E_w|=\Delta(\Delta-1)$.

Our next goal is to show that $k_{\beta}'\geq k_{\beta}$ for every $\beta$ in $[\Delta^2-1]$.

First, let $\beta$ in $F_1$ be such that $k_{\beta}=1$.
Since $F_1$ and $F_2$ are disjoint by definition, 
(d) implies that no edge in $E_{u'}$ has color $\beta$.
If $k'_{\beta}=0$, that is, no edge in $E^2_{u'}$ has color $\beta$,
then changing the color of $uv$ to $\alpha$ and the color of $uu'$ to $\beta$ 
yields a $1$-degenerate edge coloring 
with less edges colored $\Delta^2$, which is a contradiction.
Hence, $k_{\beta}'\geq 1$.

Next, let $\beta$ be a color in $F_2$ with $k_{\beta}=1$.
If $k_{\beta}'=0$, that is, no edge in $E_{u'}\cup E^2_{u'}$ has color $\beta$,
then changing the color of $uv$ to $\alpha$ and the color of $uu'$ to $\beta$ 
yields a $1$-degenerate edge coloring 
with less edges colored $\Delta^2$, which is a contradiction.
Hence, $k_{\beta}'\geq 1$.

Finally, let $\beta$ be a color in $F_2$ with $k_{\beta}=2$.
By (c), no edge in $E_{u'}$ has color $\beta$.
If $k_{\beta}'\leq 1$,
that is, there is at most one vertex in $N_{u'}$ that is incident with an edge colored $\beta$,
then changing the color of $uv$ to $\alpha$ and the color of $uu'$ to $\beta$ 
yields a $1$-degenerate edge coloring 
with less edges colored $\Delta^2$, which is a contradiction.
Hence, $k_{\beta}'\geq 2$.

Altogether, it follows that $k_{\beta}'\geq k_{\beta}$ for every $\beta$ in $[\Delta^2-1]$,
and, we obtain
\begin{eqnarray*}
\Delta(\Delta-1) 
& = & \sum\limits_{\beta\in [\Delta^2-1]}k_{\beta}
\leq \sum\limits_{\beta\in [\Delta^2-1]}k'_{\beta}
\leq |E_{u'}|+\sum\limits_{w\in N_{u'}}|E_w|
\leq (\Delta-1)+\sum\limits_{w\in N_{u'}}(\Delta-1)
\leq \Delta(\Delta-1).
\end{eqnarray*}
Equality throughout this inequality chain implies that
$k_{\beta}'=k_{\beta}$ for every $\beta\in [\Delta^2-1]$,
all edges from $E_{u'}\cup E^2_{u'}$ have a color from $[\Delta^2-1]$, and
all vertices in $N_{u'}$ have degree $\Delta$.

Now, let $v'\in N_v$.
Note that symmetric observations apply to the vertex $v'$ as to the vertex $u'$.
Let the edge $vv'$ have color $\beta$.
There is exactly one vertex $u''$ in $N_{u'}$ 
such that some edge in $E_{u''}$, say $u''u'''$, has color $\beta$.
Defining $N_{v'}$ and $E_w$ for $w\in N_{v'}$ similarly as above,
it follows, by symmetry between $u'$ and $v'$,
that there is exactly one vertex $v''$ in $N_{v'}$ 
such that some edge in $E_{v''}$, say $v''v'''$, has color $\alpha$.

If the edge $u''u'''$ is distinct from the edge $vv'$,
then changing the color of $uv$ to $\beta$ and the color of $vv'$ to $\alpha$ 
yields a $1$-degenerate edge coloring 
with less edges colored $\Delta^2$, which is a contradiction.
Hence, the edge $u''u'''$ equals $vv'$.
Since $v$ is incident with an edge colored $\Delta^2$ but $u''$ is not,
we obtain that $u''$ equal $v'$,
that is, $u'$ and $v'$ are adjacent.

\medskip

\noindent Since $u'$ and $v'$ were arbitrary vertices in $N_u$ and $N_v$, respectively,
it follows, by symmetry, that every vertex in $N_u$ is adjacent to every vertex in $N_v$,
that is, $G$ is $K_{\Delta,\Delta}$. $\Box$

\medskip

\noindent We believe that, for large values of $r$, the bound from Theorem \ref{theorem1} is far from being tight. Our next result vaguely supports this.

\begin{proposition}\label{proposition1}
If $r$ is an integer at least $2$, then no graph $G$ of maximum degree at most $\Delta$
satisfies $\chi_r'(G)=\frac{2(\Delta-1)^2}{r+1}+2(\Delta-1)+1$.
\end{proposition}
{\it Proof:} For contradiction, suppose that $G$ is a graph of maximum degree $\Delta$
that satisfies $\chi_r'(G)=K$, where $K=\frac{2(\Delta-1)^2}{r+1}+2(\Delta-1)+1$.
Similarly as in the proof of Theorem \ref{theorem2},
we consider an $r$-degenerate edge coloring of $G$ using colors in $[K]$
such that the number of edges colored $K$ is as small as possible.
Let $uv$ be an edge colored $K$.
Again using the same notation as in the proof of Theorem \ref{theorem1}
and arguing as in the proof of Theorem \ref{theorem2}, 
we obtain that 
$F_1\cup F_2=[K-1]$,
$|F_1|=2(\Delta-1)$,
$(r+1)|F_2|=2(\Delta-1)^2$,
and that the properties (a) to (g) hold.

Suppose that there is some color $\alpha$ in $F_2$ such that 
$d_u^{\alpha}$ and $d_v^{\alpha}$ are both positive.
In this case, (c) and $r\geq 2$ imply that
$\min\{ d_u^{\alpha},d_v^{\alpha}\}\leq r-1$ 
and
$\max\{ d_u^{\alpha},d_v^{\alpha}\}\leq r$,
and changing the color of $uv$ to $\alpha$
yields an $r$-degenerate edge coloring
with less edges colored $K$, which is a contradiction.
Hence, for every color $\alpha$ in $F_2$, 
we obtain, again using (c), that $(d_u^{\alpha},d_v^{\alpha})\in \{ (0,r+1),(r+1,0)\}$.

\medskip
   
\noindent Let $u'\in N_u$ and let $uu'$ have color $\alpha$.
Arguing as in the Theorem \ref{theorem2} and using the same notation as there,
it follows that every color $\beta$ in $F_1$ that appears on some edge in $E_v$
appears on at least one edge in $E^2_{u'}$.

Now, let $\beta$ be a color in $F_2$ 
such that some vertex in $N_v$ is incident with an edge colored $\beta$.
Since $d_v^{\beta}>0$ implies $d_u^{\beta}=0$ and $d_v^{\beta}=r+1$,
there are exactly $r+1$ such vertices.
If at most $r$ vertices in $N_{u'}$ are incident with an edge colored $\beta$,
then changing the color of $uv$ to $\alpha$ and the color of $uu'$ to $\beta$ 
yields an $r$-degenerate edge coloring
with less edges colored $K$, which is a contradiction.
Hence, for every such color $\beta$,
at least $r+1$ vertices in $N_{u'}$ are incident with an edge colored $\beta$,
and, since $d_u^{\beta}=0$, all these edges belong to $E^2_{u'}$.

Altogether, we obtain the contradiction
\begin{eqnarray*}
(\Delta-1)^2 
& \geq & \sum_{w\in N_{u'}}|E_w|
 \geq |E_v|+\sum_{w\in N_v}|E_w|
 = \Delta(\Delta-1),
\end{eqnarray*}
which completes the proof.
$\Box$

\section{Efficient algorithm for chordal graphs}

Let $G$ be a chordal graph.
It is well known that $G$ has a tree decomposition $(T,(X_t)_{t\in V(T)})$ such that each bag $X_t$ is a clique in $G$. By applying standard manipulations \cite{bk}, we may furthermore assume that
\begin{itemize}
\item $T$ is a rooted binary tree,
\item if $t$ is the root or a leaf of $T$, then $X_t=\emptyset$,
\item if some node $t$ of $T$ has two children $t'$ and $t''$, then $X_t=X_{t'}=X_{t''}$ ($t$ is a ``join node''), 
\item if some node $t$ of $T$ has only one child $t'$, then 

either $|X_t\setminus X_{t'}|=1$ and $|X_{t'}\setminus X_t|=0$ ($t$ is an ``introduce node'')

or $|X_t\setminus X_{t'}|=0$ and $|X_{t'}\setminus X_t|=1$ ($t$ is a ``forget node''), and
\item given $G$, the decomposition $(T,(X_t)_{t\in V(T)})$ can be constructed in polynomial time, 
in particular, $n(T)$ is polynomially bounded in terms of $n(G)$.
\end{itemize}
For every note $t$ of $T$, let $T_t$ denote the subtree of $T$ rooted in $t$ that contains $t$ and all its descendants. Let $G_t$ be the subgraph of $G$ induced by $\bigcup\limits_{s\in V(T_t)}X_s$.

We design a dynamic programming procedure calculating $\nu_r(G)$ 
for a fixed positive integer $r$.
Therefore, for every node $t$ of $T$, let ${\cal R}_t$ be the set of all triples $(S,N,k)$ such that
\begin{enumerate}[(i)]
\item $N\subseteq S\subseteq X_t$ and 
\item there is a matching $M\subseteq E(G_t)\setminus {X_t\choose 2}$ such that
\begin{enumerate}[(a)]
\item $k=|M|$,
\item $N=V(M)\cap X_t$, and 
\item $G[V(M)\cup S]$ is $r$-degenerate.
\end{enumerate}
\end{enumerate}
Note that the matching $M$ satisfying (a), (b), and (c) 
may not be uniquely determined by $(S,N,k)$.
We call every such matching {\it suitable} for $(S,N,k)$, 
and denote one (arbitrary yet specific) suitable matching by $M_t(S,N,k)$.
Intuitively, the vertices in $S$ correspond to those vertices of $X_t$ 
that can be incident with edges $e$ of some $r$-degenerate matching 
of the entire graph $G$ 
containing a suitable matching such that
either $e$ has both endpoints in $X_t$
or $e$ has one endpoint in $X_t$ and the other endpoint in $V(G)\setminus V(G_t)$.
Note that, since $X_t$ is a clique, we have $|S|\leq r+1$ for every $(S,N,k)\in {\cal R}_t$,
which implies that $|{\cal R}_t|$ is polynomially bounded in terms of $n(G)$.
Furthermore, if $t$ is the root of $T$, then 
$G_t=G$, 
all triples in ${\cal R}_t$ have the form $(\emptyset,\emptyset,k)$, and, 
by the definition of ${\cal R}_t$,
\begin{eqnarray}\label{e4}
\nu_r(G) &=& \max\left\{ k: (\emptyset,\emptyset,k)\in {\cal R}_t\right\}.
\end{eqnarray}
The following lemma contains the relevant recursions.

\begin{lemma}\label{lemma1}
Let $G$, $(T,(X_t)_{t\in V(T)})$, and $({\cal R}_t)_{t\in V(T)}$ be as above.
\begin{enumerate}[(a)]
\item If $t$ is a leaf of $T$, then ${\cal R}_t=\{ (\emptyset,\emptyset,0)\}$. 
\item If $t$ is an introduce node, $t'$ is the child of $t$, and $\{ x\}=X_t\setminus X_{t'}$, then
$(S,N,k)\in {\cal R}_t$ if and only if 
\begin{itemize}
\item either
$(S,N,k)\in {\cal R}_{t'}$
\item or 
$(S,N,k)=(S'\cup \{ x\},N,k)$ for some $(S',N,k)\in {\cal R}_{t'}$ with $|S'|\leq r$.
\end{itemize}
\item If $t$ is a forget node, $t'$ is the child of $t$, and $\{ x\}=X_{t'}\setminus X_t$, then
$(S,N,k)\in {\cal R}_t$ if and only if 
\begin{itemize}
\item either $(S,N,k)\in {\cal R}_{t'}$ and $x\not\in S$,
\item or $(S,N,k)=(S'\setminus \{ x\},N'\cup \{ y\},k'+1)$ for some $(S',N',k')\in {\cal R}_{t'}$ 
with $x\in S'\setminus N'$ and some $y\in S'\setminus (N'\cup \{ x\})$, 
\item or $(S,N,k)=(S'\setminus \{ x\},N'\setminus \{ x\},k')$ for some $(S',N',k')\in {\cal R}_{t'}$ 
with $x\in N'$.
\end{itemize}
\item If $t$ is a join node, and $t'$ and $t''$ are the children of $t$, then
$(S,N,k)\in {\cal R}_t$ if and only if 
$(S,N,k)=(S,N'\cup N'',k'+k'')$
for some $(S,N',k')\in {\cal R}_{t'}$ and $(S,N'',k'')\in {\cal R}_{t''}$ with $N'\cap N''=\emptyset$.
\end{enumerate}
\end{lemma}
{\it Proof:} (a) This follows immediately from the definition of ${\cal R}_t$.

\medskip

\noindent (b) Note that $N_{G_t}(x)=X_{t'}$,
that is, $x$ has no neighbor in $V(G_{t'})\setminus X_{t'}$.

If either $(S,N,k)\in {\cal R}_{t'}$ 
or $(S,N,k)=(S'\cup \{ x\},N,k)$ for some $(S',N,k)\in {\cal R}_{t'}$ with $|S'|\leq r$,
then the definition of ${\cal R}_t$ easily implies that $(S,N,k)\in {\cal R}_t$.
Note, in particular, that in the second case, the vertex $x$ has degree $|S'|\leq r$
in the subgraph of $G$ induced by $V(M_{t'}(S',N,k))\cup S'\cup \{ x\}$,
which ensures the degeneracy conditions.

Conversely, let $(S,N,k)\in {\cal R}_t$.
If $x\not\in S$, then, by the definition of ${\cal R}_t$, we obtain $(S,N,k)\in {\cal R}_{t'}$.
If $x\in S$, then, since $X_t$ is a clique, the set $S'=S\setminus \{ x\}$ has order at most $r$,
and, since all neighbors of $x$ belong to $X_t$, the vertex $x$ does not belong to $N$,
which implies that $(S',N,k)\in {\cal R}_{t'}$.

\medskip

\noindent (c) Note that $G_t=G_{t'}$, and that $N_G(x)\subseteq V(G_{t'})$.

If either $(S,N,k)\in {\cal R}_{t'}$ and $x\not\in S$,
or $(S,N,k)=(S'\setminus \{ x\},N'\cup \{ y\},k'+1)$ for some $(S',N',k')\in {\cal R}_{t'}$ 
with $x\in S'\setminus N'$ and some $y\in S'\setminus (N'\cup \{ x\})$, 
or $(S,N,k)=(S'\setminus \{ x\},N'\setminus \{ x\},k')$ for some $(S',N',k')\in {\cal R}_{t'}$ 
with $x\in N'$,
then the definition of ${\cal R}_t$ easily implies that $(S,N,k)\in {\cal R}_t$.
In the first case, this is immediate.
In the second case, since $x\in S'\setminus N'$ has no neighbor in $G$ outside of $V(G_{t'})$,
any suitable matching contains no edge incident with $x$
but $x$ corresponds to a vertex that can eventually be matched to some vertex $y$ in $S'\setminus (N'\cup \{ x\})$.
Since $x$ is adjacent to all vertices in $S'\setminus (N'\cup \{ x\})$, 
we add to ${\cal R}_t$ all triples corresponding to the possible choices of $y$, and increase $k'$ by $1$ because of the edge $xy$
that lies between $V(G_t)\setminus X_t$ and $X_t$.
Similarly, in the third case, the vertex $x$ is incident with an edge in $M_{t'}(S',N',k')$ whose other endpoint lies in $V(G_{t'})\setminus X_{t'}$, and removing $x$ from $X_{t'}$, it has to be removed from $S'$ and $N'$ as well while the size $k'$ of the matching $M_{t'}(S',N',k')$ does not change.

Conversely, let $(S,N,k)\in {\cal R}_t$.
Let $M=M_t(S,N,k)$.
If $x\not\in V(M)$, then $(S,N,k)\in {\cal R}_{t'}$.
If $xy\in M$ with $y\in N$, then $(S\cup \{ x\},N\setminus \{ y\},k-1)\in {\cal R}_{t'}$.
Finally, if $xy\in M$ with $y\not\in N$, 
then $(S\cup \{ x\},N\cup \{ x\},k)\in {\cal R}_{t'}$.

\medskip

\noindent (d) Note that $G_t=G_{t'}\cap G_{t''}$ and that $X_t=V(G_{t'})\cap V(G_{t''})$.

First, let $(S,N',k')\in {\cal R}_{t'}$ and let $(S,N'',k'')\in {\cal R}_{t''}$ with $N'\cap N''=\emptyset$.
Let $M=M'\cup M''$, where $M'=M_{t'}(S,N',k')$ and $M''=M_{t''}(S,N'',k'')$.
Since $N'$ and $N''$ are disjoint, 
$M$ is a matching with $M\subseteq E(G_t)\setminus {X_t\choose 2}$,
$|M|=|M'|+|M''|=k'+k''$, and
$V(M)\cap X_t=(V(M')\cup V(M''))\cap X_t=N'\cup N''$.

Let $u_1,\ldots,u_{n'}$ be a linear order of the vertices in $V(M')\cup S$ such that 
$u_1,\ldots,u_{n'-|S|}$ contains the $n'-|S|$ vertices in $V(M')\setminus N$
in an order of non-increasing depth of the corresponding forget nodes.
More precisely, 
if $1\leq i<j\leq n'-|S|$,
$t_i$ is the forget node of $u_i$, 
meaning that $u_i$ belongs $X_{t_i'}$ where $t_i'$ is the child of $t_i$ but $u_i$ no longer belongs to $X_{t_i}$, and 
$t_j$ is the forget node of $u_j$, 
then the depth of $t_i$ within $T$ is at least the depth of $t_j$.
Note that, since $(T,(X_t)_{t\in V(T)})$ is a tree decomposition, 
the forget nodes of the vertices of $G$ are uniquely determined.

Now, if $1\leq i\leq n'-|S|$, then the neighborhood of $u_i$ in the graph $G[\{ u_i,\ldots,u_{n'}\}]$
is completely contained in $X_{t_i'}$, 
because, for all vertices $u_j$ of $G_{t_i'}$ distinct from $u_i$ that belong to $V(M')\cup S$,
the forget node of $u_j$ has strictly larger depth than the forget node $t_i$ of $u_i$.
Since $X_{t_i'}$ is a clique, and $G[V(M')\cup S]$ is $r$-degenerate,
this implies that the degree of $u_i$ in the graph $G[\{ u_i,\ldots,u_{n'}\}]$ is at most $|S|-1\leq r$.
Furthermore, since $|S|\leq r+1$, for $n'-|S|+1\leq i\leq n'$,
also the degree of $u_i$ in the graph $G[\{ u_i,\ldots,u_{n'}\}]$ is at most $r$.
Altogether, it follows that $u_1,\ldots,u_{n'}$ is an $r$-degenerate order of $G[V(M')\cup S]$.
If the $r$-degenerate order $v_1,\ldots,v_{n''}$ of $G[V(M'')\cup S]$ is defined analogously, then $u_1,\ldots,u_{n'-|S|},v_1,\ldots,v_{n''}$
is an $r$-degenerate order of $G[V(M)\cup S]$,
which implies $(S,N'\cup N'',k'+k'')\in {\cal R}_t$.

Conversely, let $(S,N,k)\in {\cal R}_t$.
Let $M=M_t(S,N,k)$,
$M'=M\cap E(G_{t'})$, 
$M''=M\cap E(G_{t''})$,
$N'=V(M')\cap X_t$, and 
$N''=V(M'')\cap X_t$.
Since $M'$ and $M''$ are disjoint, 
we obtain that $k=|M'|+|M''|$ and that also the sets $N'$ and $N''$ are disjoint.
Furthermore, since the graph $G[V(M)\cup S]$ is $r$-degenerate,
also its two induced subgraphs $G[V(M')\cup S]$ and $G[V(M'')\cup S]$ are $r$-degenerate.
It follows that $(S,N',|M'|)\in {\cal R}_{t'}$ and $(S,N'',|M''|)\in {\cal R}_{t''}$.
$\Box$

\begin{theorem}\label{theorem3}
For a fixed positive integer $r$, and a given chordal graph $G$,
the maximum size of an $r$-degenerate matching can be determined in polynomial time.
\end{theorem}
{\it Proof:} By Lemma \ref{lemma1}, it follows that 
the tree decomposition $(T,(X_t)_{t\in V(G)})$ as well as the sets ${\cal R}_t$
can all be determined in polynomial time processing $T$ in a bottom-up fashion.
Furthermore, (\ref{e4}) allows to extract $\nu_r(G)$ from the set ${\cal R}_t$ of the root $t$ of $T$. 
$\Box$

\medskip

\noindent It is easy to extend the above dynamic programming approach in such a way that it also determines a maximum $r$-degenerate matching of the given graph.
Furthermore, given weights on the edges, also a maximum weight $r$-degenerate matching can be determined efficiently by replacing the cardinality $k$ within the triples $(S,N,k)$ by the weights of suitable matchings. In order to maintain the important property that the sets ${\cal R}_t$ only contains polynomially many elements, one can prune ${\cal R}_t$ maintaining only those triples $(S,N,k)$ that maximize the weight $k$ for given choices of $S$ and $N$.

Since no efficient algorithm to determine a maximum $r$-degenerate induced subgraph,
and, in particular for $r=1$, a maximum induced forest of a given chordal graph 
seems to have been published (cf. comments in \cite{kmt}),
we want to point out that modifying the above approach easily allows to obtain such an algorithm.

\section{Conclusion}

The problem to determine the acyclic matching number of a given graph $G$ 
is equivalent to the problem to determine an induced forest $T$ of $G$ 
whose matching number $\nu(T)$ is largest possible.
This observation shows that $\nu_1(G)$ is somewhat related 
to the problem to determine a largest induced forest of a given graph.
The latter problem is dual to the {\it feedback vertex set} problem,
and has received a lot of attention \cite{aks,amt,bb,gr}.
In particular, the classes of graphs for which a largest induced forest can be found efficiently 
\cite{kmt,l,lc,lt} are good candidates for classes of graphs for which the acyclic matching number might be tractable.
Note that, for a connected graph $G$ of maximum degree at most $\Delta$,
the value of $\nu_{\Delta-1}(G)$ can be determine efficiently.
In fact, if $G$ has no perfect matching, then it equals $\nu(G)$,
otherwise, it equals $\nu(G)-1$.

Further upper bounds on the $r$-degenerate chromatic index and lower bounds on the $r$-degenerate matching number for general as well as for restricted graphs seem to deserve additional research.

\end{document}